\documentclass[11pt]{article}
\usepackage[ruled]{algorithm2e}
\usepackage[top=1in, bottom=1in, left=1in, right=1in]{geometry}
\usepackage{amsmath, amsthm, amssymb}
\usepackage{color}
\usepackage{graphicx}
\usepackage{caption}
\usepackage{subcaption}
\usepackage{epstopdf}
\usepackage{overpic}
\usepackage{rotating}

\newcommand{\bXone}{\mathbf{X}}
\newcommand{\bXtwo}{\mathbf{X'}}
\newcommand{\bx}{\mathbf{x}}
\newcommand{\bb}{\mathbf{b}}
\newcommand{\bL}{\mathbf{L}}
\newcommand{\bW}{\mathbf{W}}
\newcommand{\bSigma}{\mathbf{\Sigma}}
\newcommand{\bS}{\mathbf{S}}
\newcommand{\bV}{\mathbf{V}}
\newcommand{\bPhi}{\mathbf{\Phi}}
\newcommand{\bU}{\mathbf{U}}
\newcommand{\bX}{\mathbf{X}}
\newcommand{\bA}{\mathbf{A}}
\newcommand{\bAtilde}{\mathbf{\tilde{A}}}
\newcommand{\bXS}{\bS}
\newcommand{\bXL}{\bL}
\newcommand{\balpha}{\mathbf{\alpha}}
\newcommand{\bLambda}{\mathbf{\Lambda}}

\title{\vspace{-.2in}Streaming GPU Singular Value and Dynamic Mode Decompositions}
\author{Seth D. Pendergrass \\
        \small{University of Washington} \\
        \small{Computer Engineering}
        \and
        J. Nathan Kutz \\
        \small{University of Washington} \\
        \small{Applied Mathematics}
        \and
        Steven L. Brunton \thanks{Corresponding author: +1 (609)-921-6415, sbrunton@uw.edu} \\
        \small{University of Washington} \\
        \small{Mechanical Engineering}}
\date{}

\begin{document}

\maketitle

\begin{abstract}
This work develops a parallelized algorithm to compute the dynamic mode decomposition (DMD) on a graphics processing unit using the streaming method of snapshots singular value decomposition.  
This allows the algorithm to operate efficiently on streaming data by avoiding redundant inner-products as new data becomes available.  
In addition, it is possible to leverage the native compressed format of many data streams, such as HD video and computational physics codes that are represented sparsely in the Fourier domain, to massively reduce data transfer from CPU to GPU and to enable sparse matrix multiplications.
Taken together, these algorithms facilitate real-time streaming DMD on high-dimensional data streams.
We demonstrate the proposed method on numerous high-dimensional data sets ranging from video background modeling to scientific computing applications, where DMD is becoming a mainstay algorithm.
The computational framework is developed as an open-source library written in C++ with CUDA, and the algorithms may be generalized to include other DMD advances, such as compressed sensing DMD, multi resolution DMD, or DMD with control.  

\end{abstract}

\section{Introduction}
Dynamic mode decomposition (DMD) was first introduced by Schmid in the fluids community~\cite{schmid2010dynamic} as a data-driven method to decompose complex fluid systems into spatiotemporal coherent structures, where each mode is associated with a particular frequency and rate of growth or decay.  
DMD has since been rigorously connected to nonlinear dynamical systems via Koopman operator theory~\cite{Rowley2009jfm,tu2013dynamic}, which provides an alternative infinite-dimensional linear representation of nonlinear dynamical systems~\cite{Koopman1931pnas,Mezic2005nd,Mezic2013arfm}.  
DMD may also be thought of as an algorithm~\cite{tu2013dynamic}, which yields a fundamental matrix decomposition, combining many beneficial features of principal components analysis (PCA) or proper orthogonal decomposition (POD) and the fast Fourier transform (FFT).  
As such, DMD has gained significant attention in a wide variety of fields~\cite{Kutz2016book}, including fluid dynamics~\cite{schmid2011applications,Noack2015arxiv,Brunton2015amr,Priebe2016jfm}; neuroscience~\cite{brunton2016extracting}; robotics~\cite{berger2014dynamic}; epidemiology~\cite{proctor2015discovering}; and video processing~\cite{grosek2014dynamic}.  
Despite the growing success of DMD, the underlying algorithm is based on an expensive singular value decomposition (SVD) on high-dimensional data.  
Moreover, in many applications, such as video processing and high-performance computations of transient physical processes, a windowed DMD computation must be performed repeatedly for \emph{streaming} data. 
The focus of this paper is to develop a new \emph{streaming} DMD algorithm, designed to eliminate redundant computations when repeatedly performing the DMD on a sequence of data.  

Many algorithms have been proposed to increase the speed of the SVD and DMD algorithms.  
Sayadi and Schmid~\cite{sayadi2016parallel} proposed using a parallel QR factorization as the basis for a parallel SVD on tall-skinny matrices, as are common in scientific computing and video processing. 
Hemati et al.~\cite{hemati2014dynamic} developed a batch-process and POD compressed version of the DMD, in order to accommodate large data streams. 
Brand~\cite{brand2002incremental} created the incremental SVD, a method for updating an SVD to adjust for new data.  
Brunton et al.~\cite{Brunton2015jcd} and Erichson et al.~\cite{Erichson2016cviu,erichson2015compressed} used random compression in order to reduce the size of the matrix DMD is performed on.  
In Tu et al.~\cite{tu2013dynamic} it was shown that the computational bottleneck in the DMD, when computing the singular value decomposition using the method of snapshots~\cite{sirovich1987turbulence}, is the calculation of the inner product matrix on high-dimensional data.  
They further note that when computing DMD on a sequential times-series, many redundant inner products may be avoided from one timestep to the next.   
Thus by copying these shared elements rather than recalculating them, a massive speed-up may be realized.  
The present work synthesizes and builds on many of these ideas, providing an accelerated DMD computation using a streaming method of snapshots SVD, parallelized on a GPU, and designed to work directly on natively compressed representations of the data, such as JPEG image streams.

One notable application of DMD is for the separation of \textit{background} and \textit{foreground} information from a video or sequence of images~\cite{grosek2014dynamic,erichson2015compressed}.  
In video applications, foreground/background modeling is a computationally expensive task, which only becomes more challenging with increased resolution~\cite{Bouwmans2014csr,Bouwmans2014cviu,Bouwmans2016CRC}. 
Cand{\`e}s et al.~\cite{candes2011robust} framed the problem of background subtraction as a separation of the input matrix into its sparse (foreground) and low-rank (background) components, using robust principle component analysis (RPCA).
However, RPCA is expensive, as it continues to iterate until convergence on a final result, performing a singular value decomposition on each iteration.
In contrast, DMD requires only one SVD, making it more efficient than RPCA~\cite{grosek2014dynamic} for the same task.  
Although video processing is not the primary application of DMD, it provides a challenging and intuitive set of benchmark problems to test our methods.  

\subsection{Contributions}

In this paper, we develop a streaming DMD, designed to reuse computations when processing sequential inputs.
The core of this algorithm is the streaming SVD based on the method of snapshots, which we compare to a standard SVD algorithm, demonstrating considerable speed up with negligible loss in accuracy.  
We also demonstrate a new, efficient way to calculate DMD mode amplitudes on POD coefficients, as opposed to the traditional high-dimensional least-squares fit.
Additionally, we implement both CPU and GPU versions of streaming DMD and show that these algorithms are well suited to parallel processing.
We compare the GPU implementation of the streaming DMD against a non-streaming CPU implementation, with  negligible difference in outcome.  
{Further, we design this architecture to work with the native compressed format of many data streams, including Fourier compressed image streams and the output from computational codes in the Fourier domain, to reduce data transfer from CPU to GPU and leverage sparse matrix multiplications in the streaming SVD.}  
Many of the innovations developed for streaming, GPU, compressed DMD are also equally valid for the SVD, and may have significant impact on scientific computing. 
The C++ package for the streaming DMD and SVD algorithms is available under an open-source license.  

This paper is organized as follows:  First, we review background material, including the method of snapshots SVD and the DMD in Sec.~\ref{sec:background}.
We also discuss the motivation for graphics processing unit (GPU) acceleration for our algorithms.
Next, in Sec.~\ref{sec:alg}, we explain our core innovations, including the streaming SVD/DMD, fast computation of DMD mode amplitudes, GPU acceleration, {and leveraging compressed data formats.} 
In Sec.~\ref{sec:results} we show the significant performance improvements made by our streaming algorithms, and analyze their error against the standard DMD algorithm.
Lastly, in Sec.~\ref{sec:discussion}, we summarize our findings and conclude with a  discussion on applications and future work.    

\section{Background} \label{sec:background}

In order to develop our streaming DMD algorithm, we first provide an overview of the standard DMD, the method of snapshots SVD and general purpose GPU computing.
The backbone of our streaming versions of the SVD and DMD is the method of snapshots SVD.  

In all of the analysis that follows, we consider a matrix of data snapshots $\bX\in\mathbb{R}^{n\times m}$, 
\begin{equation}
\bX = \begin{bmatrix} \vline & \vline & & \vline\\ 
\bx_1 & \bx_2 & \cdots & \bx_m\\
 \vline & \vline & & \vline   \end{bmatrix}
\end{equation}
where $n$ is the number of measurements and $m$ is the number of temporal snapshots.  For example, if the columns of $\bX$ represent image frames in a movie, then $n$ is the number of pixels per frame and $m$ is the number of frames in the movie.  Similarly, we may consider a time-series of an evolving spatial field from a numerical simulation of a partial differential equation.  

\subsection{Method of Snapshots Singular Value Decomposition (SVD)}
\begin{figure}
\begin{center}
    \begin{subfigure}{.45\textwidth}
        \includegraphics[width=\textwidth]{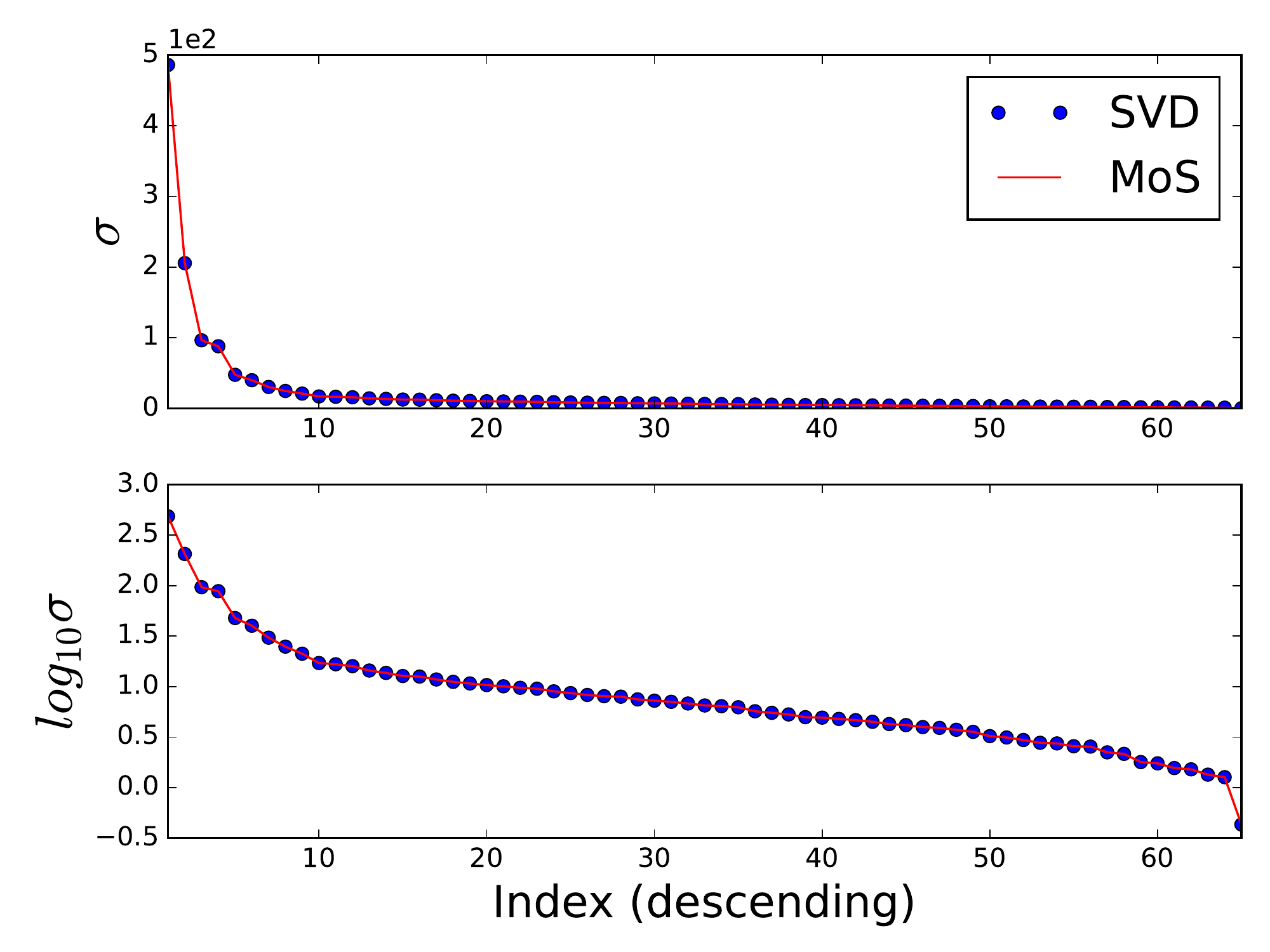}
        \caption{}
        \label{fig:sigma}
    \end{subfigure} 
    \begin{subfigure}{.45\textwidth}
        \includegraphics[width=\textwidth]{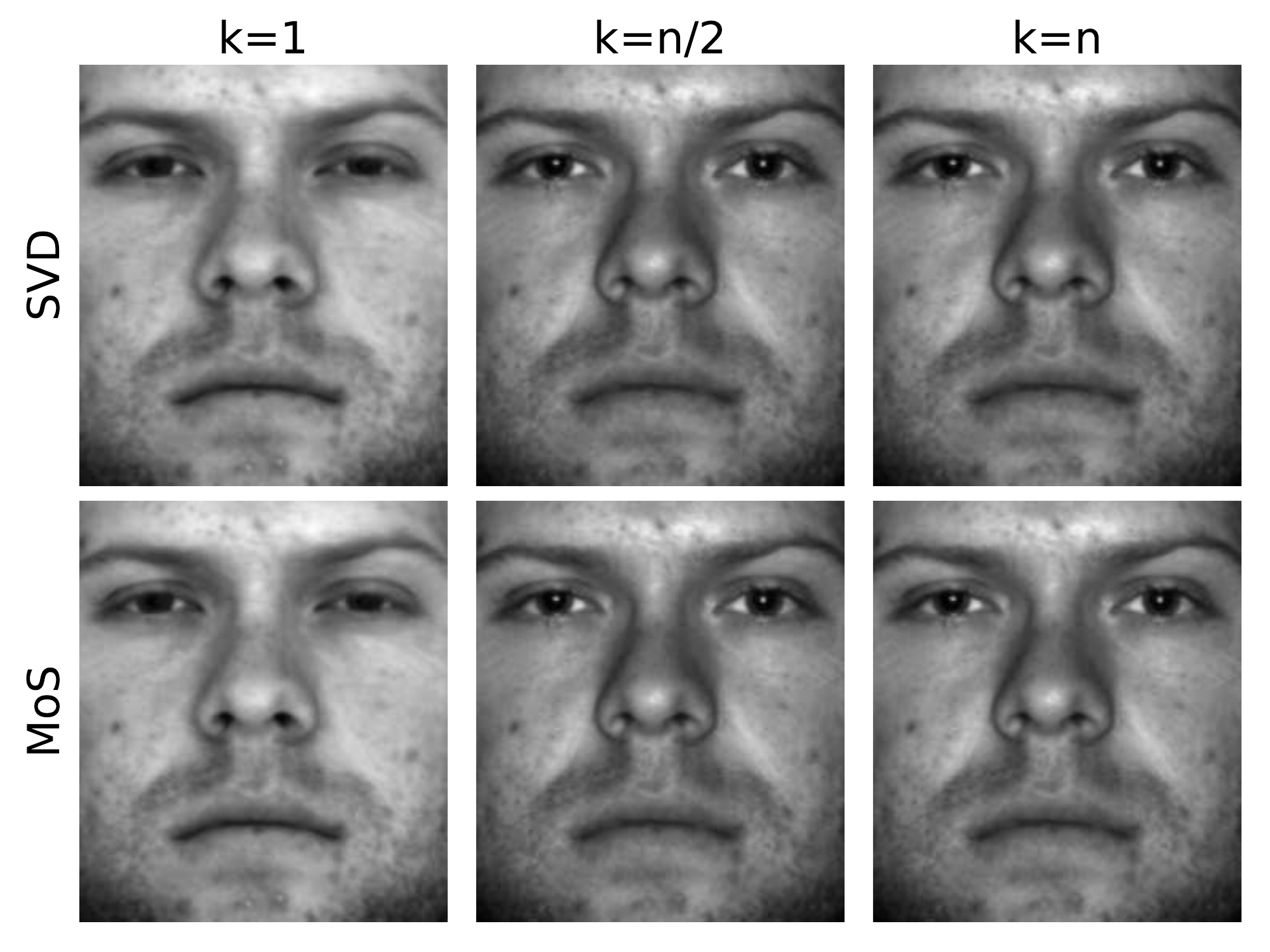}
        \caption{}
        \label{fig:faces}
    \end{subfigure}
    \vspace{-.1in}
    \caption{Comparisons between the method of snapshots SVD and the standard SVD. Relative errors between MoS and SVD for the three SVD approximations in (b) are {[2.12, 2.06, 2.06]e-7}.
             Calculating the method of snapshots SVD took .022s, while the standard SVD took .119s, a speed-up of over 5x.
             This comparison was made on the first Yale face sequence~\cite{belhumeur1997eigenfaces}.
             (a) is a comparison of singular values from the standard and method of snapshots SVDs.
             Blue circles are from the standard SVD, while the red line from the method of snapshots.
             (b) is a comparison of the standard and method of snapshots SVDs for image reconstruction with one, half or all singular values and vectors.}
             \vspace{-.2in}
    \label{fig:svd}
    \end{center}
\end{figure}

The method of snapshots is an alternative way to calculate the singular value decomposition of a matrix $\bX$, 
\begin{equation}
\bX = \bU\bSigma\bV^*,
\end{equation}
developed for matrices where one dimension is much larger than the other.  
This method was originally developed for data from fluid dynamics, in which the target matrices are significantly taller than they are wide~\cite{sirovich1987turbulence}, i.e. $n\gg m$.  
In these applications, it is observed that the nonzero eigenvalues of $\bX^*\bX$ are the same as those of $\bX\bX^*$, although the first matrix is size $m\times m$ while the second matrix is size $n\times n$.  
It is computationally more efficient to compute the eigendecomposition of the smaller matrix $\bX^*\bX$ and then use this information to reconstruct the left and right singular vectors of $\bX$.   
This allows for significant reductions in computation time, although with a potential reduction in accuracy.
The method of snapshots is summarized as follows:
\begin{enumerate}
    \item Multiply $\bX$ by its transpose, in whichever order creates the smallest output.
    We assume $\bX$ is a tall-skinny matrix (i.e., $n\gg m$).
    Then find the eigendecomposition:
    \begin{align}
        \bX^*\bX\bV=\bV\bLambda
    \end{align}
    where $\bLambda$ are the eigenvalues and $\bV$ the eigenvectors of $\bX^* \bX$.
    The non-negative square roots of $\bLambda$ are the singular values $\bSigma$ of the original matrix $\bX$.

    \item The left singular vectors $\bU$ are calculated as follows:
    \begin{align}
        \bU=\bX\bV\bSigma^{-1}.
    \end{align}
\end{enumerate}
This creates an ``economy" SVD, where $\bU\in\mathbb{R}^{n\times m}$ is the same dimension as $\bX$, and $\bSigma\in\mathbb{R}^{m\times m}$ and $\bV\in\mathbb{R}^{m\times m}$ are both small square matrices.
Figure \ref{fig:svd}(a) shows the singular values calculated with both the standard SVD and the method of snapshots, performed on the Yale faces dataset~\cite{belhumeur1997eigenfaces}.  
The singular values from both methods agree closely; thus it is only when extreme accuracy is needed that the standard method for calculating the SVD should be used.
We believe most users of the SVD would benefit by using the method of snapshots, due to its significantly better performance.  
The method of snapshots is a standard technique in the fluid dynamics community due to the high aspect ratio of the data matrix.  

Figure \ref{fig:svd}(b) compares the reconstruction of the Yale faces~\cite{belhumeur1997eigenfaces} between the standard SVD and method of snapshots SVD, as well as the absolute difference between the two.
This further demonstrates how close the method of snapshots is to the standard SVD, irregardless of the number of eigenvalues and eigenvectors used to reconstruct the images.
In turn, the speed-up provided by the method of snapshots SVD can be carried over to the DMD.

\subsection{Dynamic Mode Decomposition}
The DMD arose out of the fluid dynamics community to analyze the spatio-temporal coherent structures arising from fluids data~\cite{schmid2010dynamic}.  
It quickly gained popularity as strong connections were made between DMD and Koopman spectral analysis~\cite{Rowley2009jfm,tu2013dynamic,Kutz2016book}, which provides an infinite-dimensional linear representation of nonlinear dynamical systems~\cite{Koopman1931pnas,Mezic2005nd,Mezic2013arfm}.

DMD finds the dominant eigenvalues and eigenvectors of a best-fit linear dynamical system modeling the transition of a state $\bx_k$ to the next time-step $\bx_{k+1}$; nonlinear model reduction is also possible with similar data~\cite{Brunton2016pnas}.  In particular, given a matrix $\bX$ and another matrix $\bX'$ consisting of the snapshots one time-step in the future:
\begin{equation}
\bX' =   \begin{bmatrix} \vline & \vline & & \vline\\ 
\bx_2 & \bx_3 & \cdots & \bx_{m+1}\\
 \vline & \vline & & \vline   \end{bmatrix},
\end{equation}
the DMD algorithm obtains the eigendecomposition of the best-fit linear operator $\bA$ given by
\begin{equation}
\bA = \bX'\bX^{\dagger} = \bX' \bV\bSigma^{-1}\bU^*,
\end{equation}
where $^{\dagger}$ denotes the Moore-Penrose pseudo inverse~\cite{tu2013dynamic}.  

However, since the state dimension $n$ may be quite large (on the order of a million for HD video, tens of millions for 4K video, and even larger for scientific computing applications), the matrix $\bA$ is too large to directly analyze on simple computational architectures.  
Instead, it is possible to analyze a smaller matrix $\bAtilde$ obtained via projection onto the left singular vectors in $\bU$:
\begin{align}
    \bAtilde  = \bU^* \bA \bU =  \bU^*\bXtwo\bV\bSigma^{-1}.\label{Eq:Atilde}
\end{align}
Much like the method of snapshots, the matrix $\bAtilde$ is size $m\times m$, and it has the same eigenvalues as the high-dimensional matrix $\bA$, as shown in~\cite{tu2013dynamic}.  Taking the eigendecomposition
\begin{align}
    \bAtilde\bW=\bW\bLambda
\end{align}
it is then possible to obtain eigenvectors of the original high-dimensional matrix $\bA$ via
\begin{align}
    \bPhi=\bXtwo\bV\bSigma^{-1}\bW.\label{Eq:Phi}
\end{align}
The columns of $\bPhi$ are called \emph{dynamic modes} of $\bX$ and they are spatio-temporal modes that have a single temporal signature given by the corresponding eigenvalue $\lambda$ in $\bLambda$.  

The large number of independent inner-products performed in the process of calculating the SVD and DMD make it a perfect fit for being computed on a Graphics Processing Unit (GPU), where their many cores can be leveraged.

\subsubsection{DMD for Video Background Subtraction}
Grosek and Kutz~\cite{grosek2014dynamic} show that the DMD can be effectively leveraged to compute decomposition of a video into the foreground and background components.  
This provides a similar decomposition as in the robust principle component analysis (RPCA)~\cite{candes2011robust}, but at a fraction of the cost, as RPCA involves an iterative procedure requiring dozens of SVD computations.  In this framework, the video $\bX$ is decomposed into its constituent \textit{low-rank} and \textit{sparse} components, where the \textit{low-rank} contains a low-dimensional representation of the system under observation and the \textit{sparse} the outliers, noise and/or corruption measured by the input.
This is represented as:
\begin{align}
    \bX = \bXL + \bXS,
\end{align}
where $\bXL$ is the low-rank component (background) and $\bXS$ is the sparse component (foreground).  

Because each DMD mode has a corresponding frequency given by the DMD eigenvalue $\lambda$, the discrete-time eigenvalues that are nearly equal to $1$ correspond to modes that do not change from frame to frame, i.e., the background modes.  
Thus, DMD can also be used to split the matrix $\bX$ into two components, corresponding to slowly varying modes with eigenvalues $\lambda_p\approx 1$, and those that have faster dynamics:
\begin{align}
\bX = \underbrace{\sum_{p} b_p\phi_p\lambda_p^{\mathbf{t}-1}}_{\text{Background modes}} +\underbrace{\sum_{j\neq p} b_j\phi_j\lambda_j^{\mathbf{t}-1}}_{\text{Foreground modes}},
\end{align}
where $\mathbf{t} = \begin{bmatrix} 1 & 2 & \cdots & m\end{bmatrix}$ is a vector of time indices.  
Refer to Erichson et al. for the state of the art DMD implementation of background modeling~\cite{erichson2015compressed}.

\subsection{General Purpose GPU Computing}

General purpose GPU (GPGPU) computing has proven effective for accelerating many linear algebra problems, as it is able to perform many operations in parallel.
Creating an efficient algorithm for use on a modern GPU requires a very different approach than would be used on a central processing unit (CPU).
A typical GPU is made up of a number of sub-processors, each able to run multiple threads concurrently.
This design allows a GPU to achieve a much higher throughput than a CPU~\cite{nvidiacuda2}, if the algorithm is written with the GPU in mind.
This \textit{Single-instruction, Multiple Data (SIMD)} style of code works best when there is a large amount of input data needing to be independently processed.
Matrix multiplication is a common example, however one could expect many large-scale math problems to suit the GPU architecture well.
NVIDIA~\cite{nvidiacuda2} notes that minimizing host (CPU) to device (GPU) memory transfers is key to maximizing performance.
This lends itself naturally, then, to streaming algorithms, where only the updated data need be transferred on/off the device.

\section{Proposed Streaming DMD Algorithm} \label{sec:alg}

In many applications, data is continually acquired from sensors in a streaming fashion; new data is appended as columns to the right of the matrix $\bX$, while old columns may be removed from the left of $\bX$ if necessary.  
In streaming applications, such as online video processing or windowed DMD on transient simulations, the cost of repeated DMD and SVD calculations may be prohibitively expensive.  

Here, we build a suite of complementary techniques to accelerate repeated SVD and DMD computations for streaming data.  
The core of the streaming DMD algorithm is the streaming method of snapshots SVD, whereby redundant inner product computations in $\bX^*\bX$ are reused from one timestep to the next, reducing the SVD computational complexity from $\mathcal{O}(nm^2)$ to $\mathcal{O}(nm)$.  
The streaming SVD and DMD are discussed in Secs.~\ref{Sec:Methods:SVD} and~\ref{Sec:Methods:DMD}, respectively.  
When it is necessary to compute the mode amplitudes in $\bb$, we introduce an efficient computation in Sec.~\ref{Sec:Methods:ComputeB}.  
All of the above methods are readily parallelized, and we discuss GPGPU implementation in Sec.~\ref{Sec:Methods:GPU}.  
{Once GPU parallelized algorithms have been implemented, data transfer from the CPU to GPU becomes the main computational bottleneck.  However, in many applications it is possible to leverage the native sparse representation of the data (e.g., image sequences are stored in compressed Fourier or wavelet representations) to significantly reduce data transfer and promote sparse matrix operations, further reducing the computational burden.  This is discussed in Sec.~\ref{Sec:Methods:Compression}.}

\begin{figure}
    \centering
    \vspace{-.15in}
    \includegraphics[width=.8\textwidth]{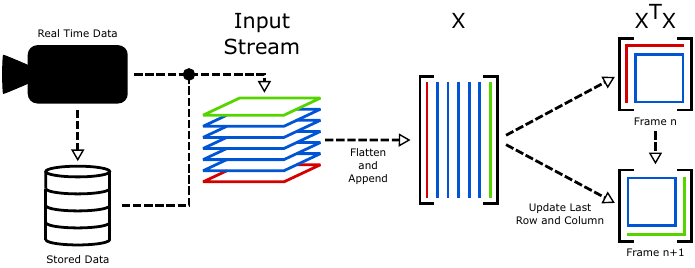}
    \vspace{-.05in}
    \caption{Overview of the streaming method of snapshots SVD, which avoids redundant inner product computations by reusing portions of the matrix $\bX^*\bX$ in subsequent time steps.}
    \vspace{-.1in}
    \label{fig:streaming}
\end{figure}

\subsection{Streaming Method of Snapshots SVD}\label{Sec:Methods:SVD}

In the streaming context, let $\bX$ be the current data matrix and $\bXtwo$ be the next matrix in the sequence.  Many of the inner products in $\bX^*\bX$, shown in blue, may be reused in $\bXtwo^*\bXtwo$:
\begin{subequations}
\begin{eqnarray}
    \bXone^*\bXone &=&
    \begin{bmatrix}
        {\color{red}\langle \bx_1,\bx_1\rangle} & {\color{red}\langle \bx_1,\bx_2\rangle} & {\color{red}\cdots} & {\color{red}\langle \bx_1,\bx_{n-1}\rangle}\\
        {\color{red}\langle \bx_2,\bx_1\rangle} & {\color{blue} \langle \bx_2,\bx_2\rangle} & {\color{blue} \cdots} & {\color{blue}\langle \bx_2,\bx_{n-1}\rangle}\\
        {\color{red}\vdots} & {\color{blue}\vdots} & {\color{blue}\ddots} & {\color{blue}\vdots} \\
        {\color{red}\langle \bx_{n-1},\bx_1\rangle} & {\color{blue}\langle \bx_{n-1},\bx_2\rangle} & {\color{blue}\cdots} & {\color{blue}\langle \bx_{n-1},\bx_{n-1}\rangle}
    \end{bmatrix}\\
    \bXtwo^*\bXtwo &=&
    \begin{bmatrix}
        {\color{blue} \langle \bx_2,\bx_2\rangle} &{\color{blue} \cdots} & {\color{blue}\langle \bx_2,\bx_{n-1}\rangle} & {\color{green}\langle \bx_2,\bx_{n+1}\rangle}\\
        {\color{blue}\vdots} & {\color{blue}\ddots} & {\color{blue}\vdots} & {\color{green}\vdots} \\
        {\color{blue}\langle \bx_{n-1},\bx_2\rangle }& {\color{blue}\cdots} & {\color{blue}\langle \bx_{n-1},\bx_{n-1}\rangle} & {\color{green}\langle \bx_{n-1},\bx_{n}\rangle}\\
        {\color{green}\langle \bx_{n},\bx_2\rangle} &  \cdots &{\color{green}\langle \bx_{n},\bx_{n-1}\rangle} & {\color{green}\langle \bx_{n},\bx_{n}\rangle}
    \end{bmatrix}.
\end{eqnarray}
    \end{subequations}
Thus, as $\bXtwo^*\bXtwo$ is symmetric, only the last row or column will need to be recalculated (shown in green).    
Removing the redundant inner product calculations reduces the computational complexity from $O(nm^2)$ to $O(nm)$.
As this is the most time-consuming part of the method of snapshots~\cite{tu2013dynamic}, a large performance gain is realized.
This streaming method of snapshots facilitates a streaming version of the DMD.  

\subsection{Streaming Dynamic Mode Decomposition}\label{Sec:Methods:DMD}

\begin{figure}
    \centering
    \includegraphics[width=.8\textwidth]{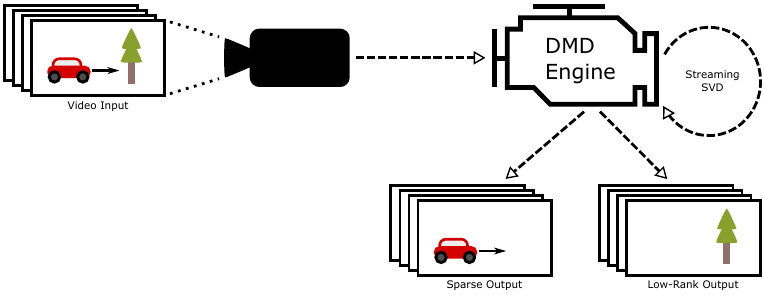}
    \caption{Diagram of the streaming DMD engine used for video separation.}
    \label{fig:engine}
\end{figure}

The streaming DMD relies on the streaming SVD in order to process data in sequence, but is also able to realize speed-ups from reusing intermediate steps from the SVD, and by only returning the last column of the sparse matrix $\bXS$ in the case of background subtraction.
Figure \ref{fig:engine} shows an outline of how the streaming DMD is set up in order to perform background subtraction. 
Background subtraction with the streaming DMD is performed by sliding the DMD forward by as many frames as the user wants to process in a given iteration.  
However, the same process would also be used if only the DMD modes and their frequencies are desired.  

\subsection{Computing Mode Amplitudes Efficiently}\label{Sec:Methods:ComputeB}

Computing the vector $\bb$ of DMD mode amplitudes has been investigated in the past~\cite{Jovanovic2014pof,Kutz2016book}.  The simplest approach involves computing a best-fit $\bb$ vector using the least-squares approximation:
\begin{align}
    \bb = \bPhi^{\dagger}\bx_1.
\end{align}
Instead, we use the following formulation directly on POD coefficients using Eqs.~\eqref{Eq:Atilde} and~\eqref{Eq:Phi}:
\begin{subequations}
\begin{align}
& \bx_1  = \bPhi \bb\\
\quad\Longrightarrow\quad & \bU \balpha_1 = \bX'\bV\bSigma^{-1}\bW \bb\\
\quad\Longrightarrow\quad & \balpha_1 = \tilde{\bA}\bW \bb\\
\quad\Longrightarrow\quad & \balpha_1 =\bW \bLambda\bb\\
\quad\Longrightarrow\quad & \bb = \left(\bW\bLambda\right)^{-1}\balpha_1,
\end{align}
\end{subequations}
where $\balpha_1$ is the vector of POD coefficients for $\bx_1$.  This is significantly more efficient than the high-dimensional least-squares algorithm.  
Additionally, only the row corresponding to the smallest absolute DMD eigenvalue need be calculated when streaming.
The benefit of this faster calculation of the DMD mode amplitudes is even more pronounced on a GPU, requiring fewer synchronizations with the device, and reducing the amount of data transfer.  

\subsection{GPU Implementation}\label{Sec:Methods:GPU}
Algorithms \ref{alg:svd}, \ref{alg:dmd} and \ref{alg:bs} show how the calculations for the SVD, DMD and background subtraction are performed.
$\bU$ is not explicitly calculated so as to reduce space and computational complexity of the DMD.
While it is possible for this to cause issues with numerical accuracy, we found the results of these algorithms to be negligibly less accurate.   
Additionally, we used single-precision to further reduce memory usage and increase performance.
Our code relies on OpenBLAS~\cite{wang2013augem} for LAPACK and BLAS functions on the CPU, and MAGMA~\cite{tdb10} for GPU LAPACK and BLAS.  
We also found that writing algorithms in MAGMA improved performance over those written in OpenCL. 

\begin{algorithm}
    \DontPrintSemicolon
    \caption{Streaming Method of Snapshots SVD}
    \label{alg:svd}
    \SetKwProg{Fn}{Function}{:}{\KwRet{sigma, v, xtx}}
    \SetKwFunction{FSVD}{SSVD}
    \Fn{\FSVD{X, first\_step}}{
        \eIf{first\_step}{
            xtx = X.T * X\;
        }{
            xtx[:-1, :-1] = xtx[1:, 1:]\;
            xtx[:, -1] = X.T * X[:, -1]\;
            xtx[-1, :] = xtx[:, -1].T\;
        }
        s, v = eig(xtx)\;
        sigma = sort(sqrt(abs(s)), 'desc')\;
    }
\end{algorithm}
\vspace{-.15in}
\begin{algorithm}
    \DontPrintSemicolon
    \caption{Streaming DMD}
    \label{alg:dmd}
    \SetKwProg{Fn}{Function}{:}{\KwRet{phi, lambda, sigma, v, w}}
    \SetKwFunction{FDMD}{SDMD}
    \Fn{\FDMD{X, first\_step}}{
        sigma, v, xtx = SSVD(X[:, :-1], first\_step)\;
        xty[:, :-1] = xtx[:, 1:]\;
        xty[:, -1] = X[:, :-1].T * X[:, -1]\;
        vsi = v * inv(sigma)\;
        atilde = vsi.T * xty * vsi\;
        lambda, w = eig(atilde)\;
        vsiw = vsi * w\;
        phi = X[:, 1:] * vsiw\;
    }
\end{algorithm}
\vspace{-.15in}
\begin{algorithm}
    \DontPrintSemicolon
    \caption{Streaming DMD Background Subtraction}
    \label{alg:bs}
    \SetKwProg{Fn}{Function}{:}{\KwRet{s}}
    \SetKwFunction{FBackSub}{SBackSub}
    \Fn{\FBackSub{X, first\_step}}{
        phi, lambda, sigma, v, w = SDMD(X, first\_step)\;
        alpha1 = sigma * v[:, 0].T\;
        wl = w * lambda\;
        b = lstsq(wl, alpha1)\;
        idx = argmin(i, abs(log(lambda[i])))\;
        \eIf{first\_step}{
            lambdaPow = pow(lambda[idx], [0:X.shape[1]])\;
            l = b[idx] * phi[:, idx] * lambdaPow\;
            s = X - abs(l)\;
        }{
            lambdaPow = pow(lambda[idx], X.shape[1])\;
            l = b[idx] * phi[:, idx] * lambdaPow\;
            s = X[:, -1] - abs(l)\;
        }
    }
\end{algorithm}

\subsection{Implementation on Sparse Data}\label{Sec:Methods:Compression}

\begin{figure}
\begin{center}
\includegraphics[width=.675\textwidth]{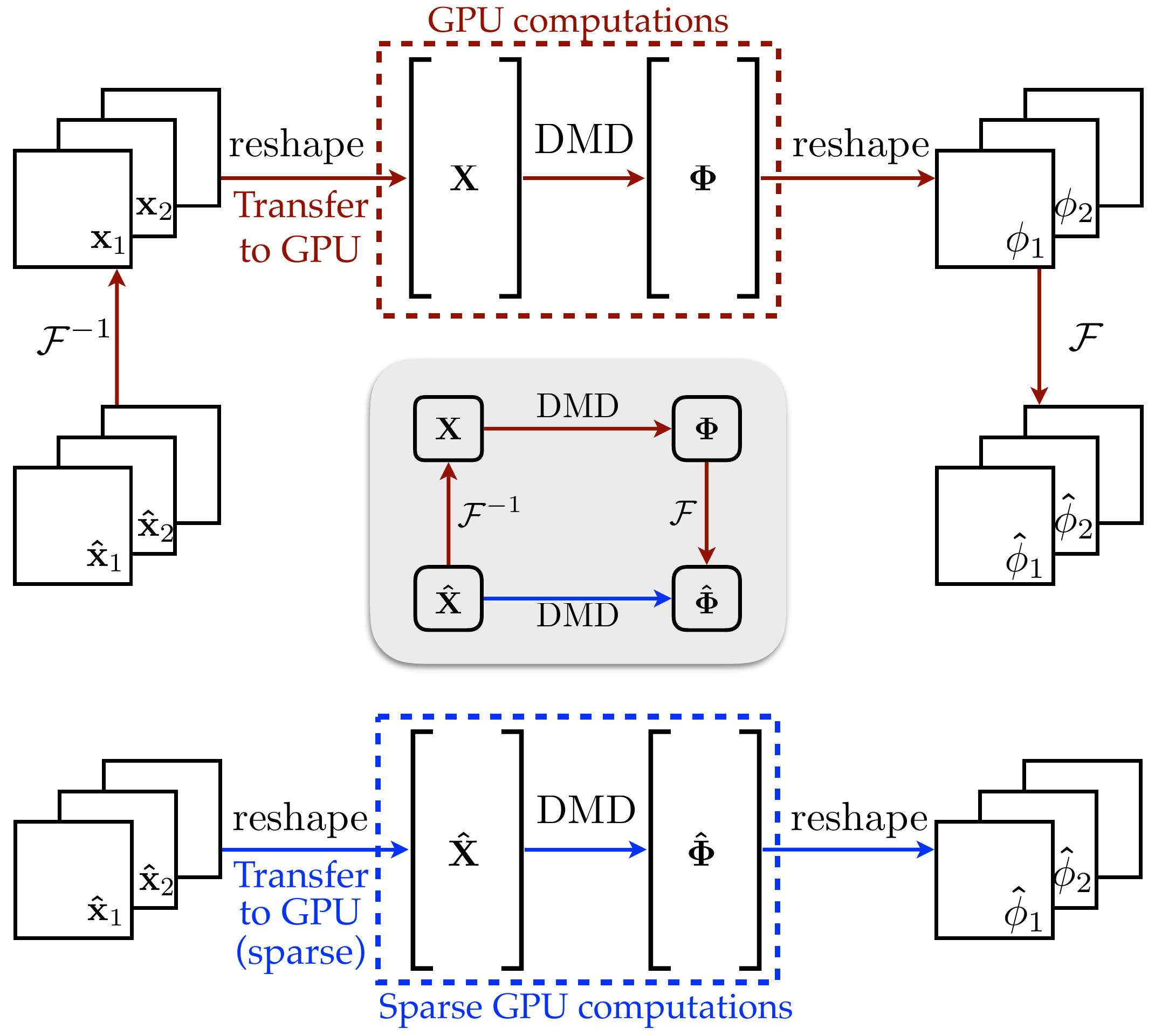}
\end{center}
\vspace{-.25in}
\caption{{Illustration of GPU accelerated DMD leveraging native compressed data formats, such as for image streams that are compressed in the Fourier domain.}}\label{Fig:CompressedDMD}
\vspace{-.15in}
\end{figure}

After parallelization on the GPU, data transfer from the CPU to the GPU and back becomes a bottleneck.  
We may naively transfer data in the ambient signal space, such as pixel space for images or a spatial domain for high performance computations.  
However, in both cases, these signals are typically stored or computed in a transformed basis, such as Fourier or wavelets.  
Moreover, these transform bases allow the data to be massively compressed, often by orders of magnitude, which would lead to a significant savings in data transfer.  
Recent work combining compressed sensing and DMD~\cite{Brunton2015jcd} showed that both the SVD and DMD  are invariant to unitary transformation, such as the fast Fourier transform (FFT).  
Thus, it is possible to directly transfer FFT compressed data to the GPU, perform DMD on the Fourier representation, and transfer the compressed DMD from the GPU back for storage or further processing.    
There is an added benefit that many of the core steps in the DMD algorithm will be performed on sparse data matrices, enabling further efficiency gains.  
This procedure is shown schematically in Fig.~\ref{Fig:CompressedDMD}.   
This is not explicitly implemented in our code, but is included because of the potentially important  role in reducing transfer from CPU to GPU in practical implementations.  
Note that compressed and randomized~\cite{Halko2011siamreview} architectures have recently been used to great advantage in scientific computing applications, for example in~\cite{Schaeffer2013pnas,mackey2014compressive}.

\section{Results} \label{sec:results}
We now present the performance and accuracy comparisons of our streaming SVD, DMD and background subtraction algorithms. 
The algorithms are demonstrated on high-resolution video data because they are publicly available and reproducible.  
However, the streaming DMD method is general to any high-dimensional data, such as data generated by high-performance computing, internet of things, LIDAR sensors, etc.  
We use the PEViD ``Walking Day Indoor 4" video~\cite{KorUHD2014} to test performance and scaling, and to subjectively compare the results of background subtraction.
The high resolution of this sequence allows us to explore the ability of our algorithm to scale with varying resolution. 
Additionally, we use the BMC ``Video 003" to make a quantitative analysis of the accuracy of our streaming GPU DMD for background subtraction in terms of standard metrics in computer vision. 

\subsection{Performance Benchmarks}

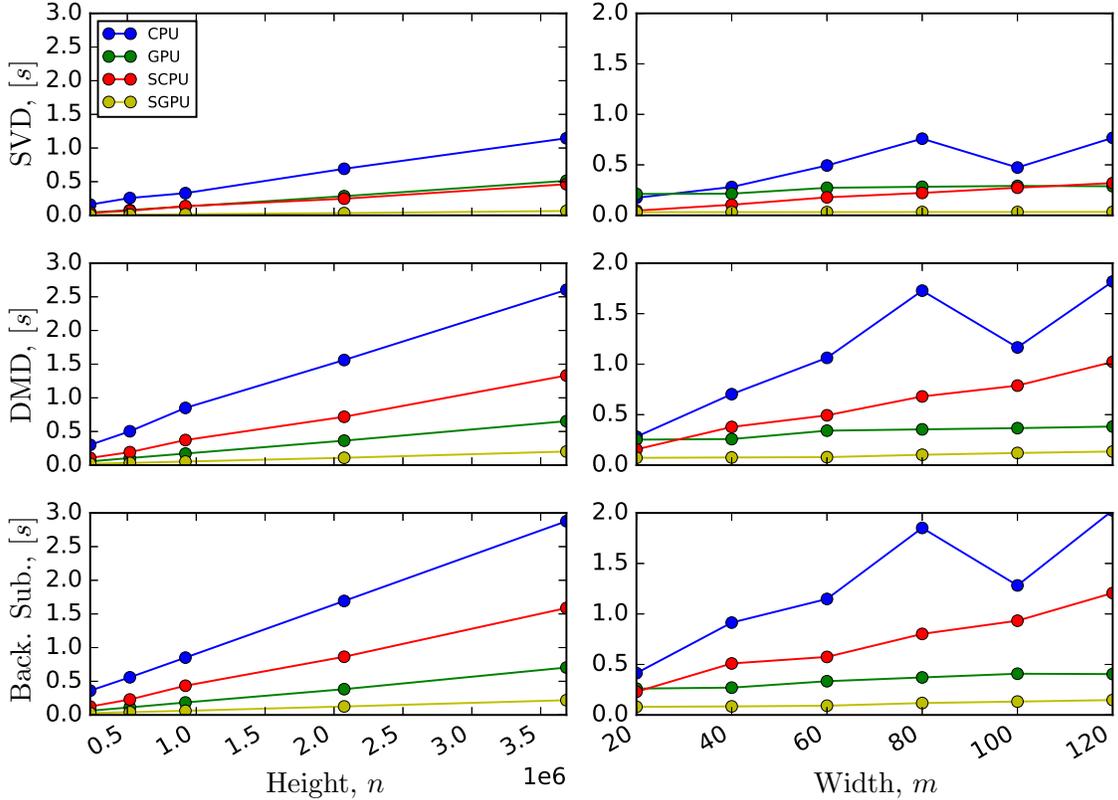
\begin{figure}
\begin{center}
	\begin{overpic}[width=.9\textwidth]{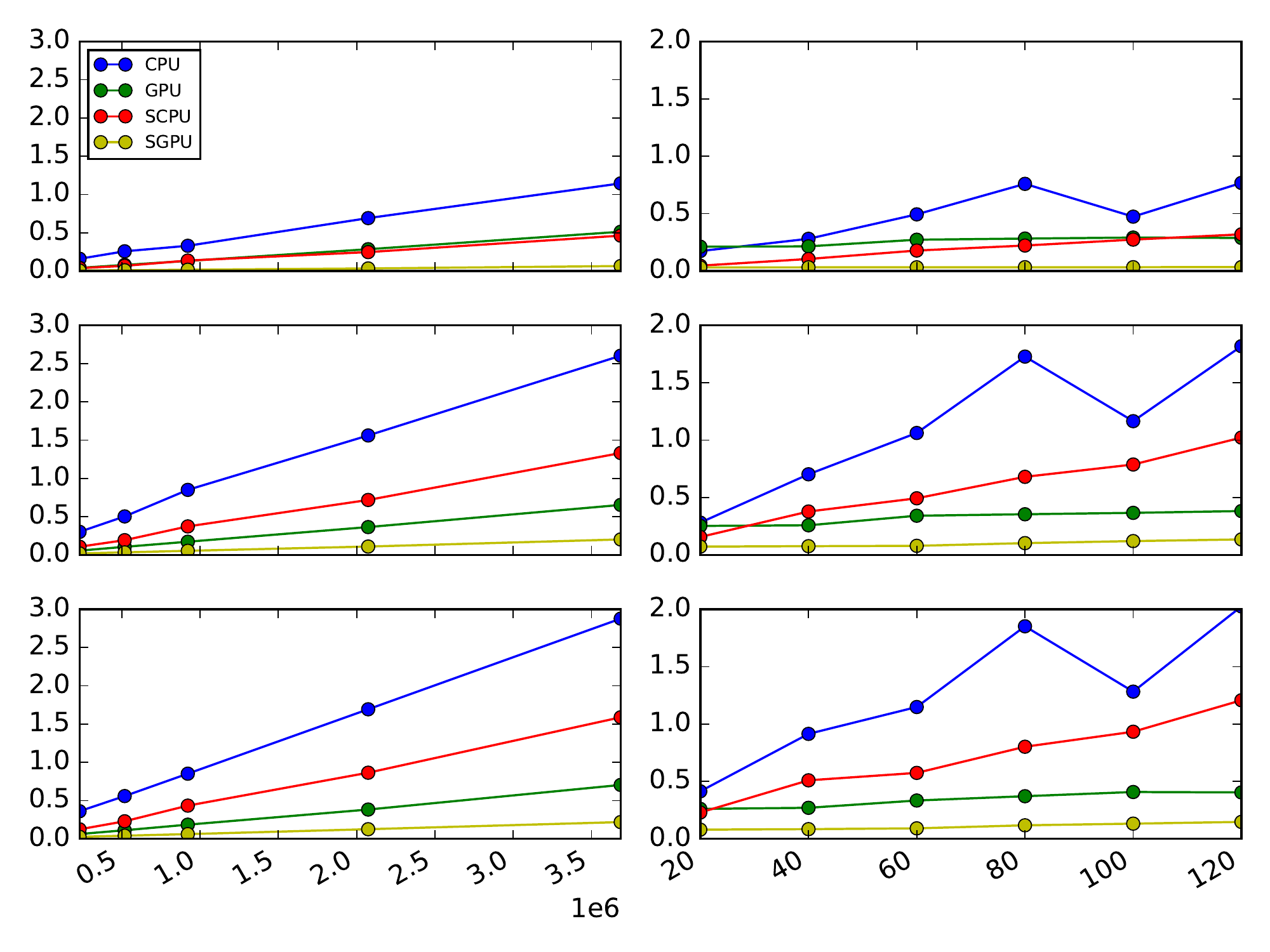}
		\put(-1,10){\begin{sideways}Back. Sub., [$s$]\end{sideways}}
		\put(-1,35){\begin{sideways}DMD, [$s$]\end{sideways}}
		\put(-1,58){\begin{sideways}SVD, [$s$]\end{sideways}}
		\put(22,2){Height, $n$}
		\put(71,2){Width, $m$}
	\end{overpic}
	\end{center}
	\vspace{-.2in}
    \caption{Comparison of CPU, streaming CPU (SCPU), GPU and streaming GPU (SGPU) versions of the SVD, DMD and DMD background subtraction.
             Tests are run with a constant width of $m=90$ (left) or input height of $n=1920 \times 1080$ (right) on PEViD ``Walking Day Indoor 4"~\cite{KorUHD2014}.}
    \label{fig:perf}
\end{figure}
We benchmark the various algorithms on the PEViD ``Walking Day Indoor 4" video~\cite{KorUHD2014}, converted to greyscale and resized to common 16:9 resolutions.\footnote{Our tests were performed on an Intel Xeon E5-2620v3 with 32GB of RAM and an NVIDIA Tesla K40, running Ubuntu 16.04.1 LTS.
          Code was compiled with gcc 5.4.0 and flags -O3 -march=native, with dependencies on OpenBLAS 0.2.19~\cite{wang2013augem} and MAGMA 2.1.~\cite{tdb10}}
We choose not to include data transfers between CPU and GPU in our benchmarks, as they don't reflect the computational differences in CPU versus GPU code; instead they are only representative of a hardware limitation of current computers.  
However, in practical implementations, we discuss the potential to significantly reduce data transfer using compressed data formats as shown in Fig.~\ref{Fig:CompressedDMD} in Sec.~\ref{Sec:Methods:Compression}.  
Instead, our tests measure the time taken to update the SVD, DMD or DMD background subtraction from steady state, where the system has already been initialized. 
For both implementations, the initial time would be equal to the elapsed time taken by their respective non-streaming versions to update.

Figure \ref{fig:perf} shows comparisons of the CPU, GPU, streaming CPU (SCPU) and streaming GPU (SGPU) implementations for the SVD, DMD and DMD background subtraction.
Our first set of benchmarks holds the width constant at $m=90$ frames (number of columns in ${\bf X}$), and varies the resolution from $n=640\times 360$ to $n=2560\times1440$ (number of rows in ${\bf X}$).
Likewise, resolution is kept constant at $n=1920\times 1080$, and width is varied from $m=20$ to $m=120$ in steps of 20 in our second benchmark set.
This test shows streaming to significantly benefit the CPU implementation, putting it on par with that of the GPU.
In real world applications, this is promising as it could reduce the need for a dedicated GPU, while still netting a large performance improvement. 
Further, the streaming GPU is significantly faster than the other three versions, and has a much smaller slope, suggesting better scaling for even larger input dimensions (i.e., resolution $n$ and number of frames $m$).
This trend is maintained in the DMD and background subtraction algorithm benchmarks as well.
Looking at the $n=2560\times 1440$ resolution by $m=90$ frames test on the SVD, the CPU implementation took approximately 1.15 seconds, while the streaming GPU took only .06 seconds, for a speed-up of nearly 20x.  
When the resolution is kept constant, we see that the scaling of both streaming algorithms is more favorable than that of the non-streaming algorithms.
This is as expected, since that the cost to update $\bX^*\bX$ is on the order of $O(n)$.
This shows that the streaming algorithm scales more favorably with regards to width than the standard SVD and DMD, as well as for the DMD for background subtraction.
Similar to the constant width tests above, in the constant height test, a speed-up of around 25x is realized for the SVD at a $m=120$ frame width, with the CPU taking .77 seconds and the streaming GPU taking .03 seconds.

\subsection{Error Analysis}
\begin{table}
	\begin{center}
		\begin{tabular}{|l|l|l|l|l|l|l|l|}
			\hline
			Dataset & Alg. & $m=20$  & $m=40$ & $m=60$ & $m=80$ & $m=100$ & $m=120$ \\
			\hline
			PEViD & SVD & .00197074 & .00176278 & .00178229 & .00164394 & .00155754 & .00240962 \\
			& DMD & .0215632 & .00878696 & .00762604 & .00906957 & .00684136 & .0293789 \\
			\hline
			BMC& SVD & 1.91137e-6 & 1.17378e-5 & 2.22748e-5 & 5.91658e-5 & .000112305 & 9.45566e-5 \\
			& DMD & 5.92535e-6 & 4.02722e-5 & 5.30704e-5 & .000191003 & 8.95717e-5 & 7.28629e-5 \\
			\hline
		\end{tabular}
	\end{center}
	\caption{Relative errors of SVD ($\bSigma$) and DMD ($\phi_{\lambda\approx0}$) streaming GPU implementation versus standard Python implementation. Measurements were made using the PEViD~\cite{KorUHD2014} dataset at $n=1920\times 1080$ height and the BMC~\cite{Vacavant2013} dataset at varying widths from $m=20$ to $m=120$.}
	\label{tab:error}
\end{table}

\begin{table}
	\begin{center}
		\begin{tabular}{|l|l|l|l|l|}
			\hline
			Sequence & Recall & Precision & F-measure & Psnr \\
			\hline
			BMC Video 003 & .713289 & .854605 & .777579 & 41.6178 \\
			\hline
		\end{tabular}
	\end{center}
	\caption{Results created using BMC Evaluation Wizard on the results of streaming GPU DMD for background subtraction on BMC ``Video 003"~\cite{Vacavant2013}. A width of 60 frames was used for streaming. The foreground mask was generated by setting all values less than .2 to 0 and any greater to 1 after subtracting the low-rank matrix from the original input.}
	\label{tab:video}
\end{table}

\begin{figure}\begin{center}
\begin{overpic}[height=.275\textwidth]{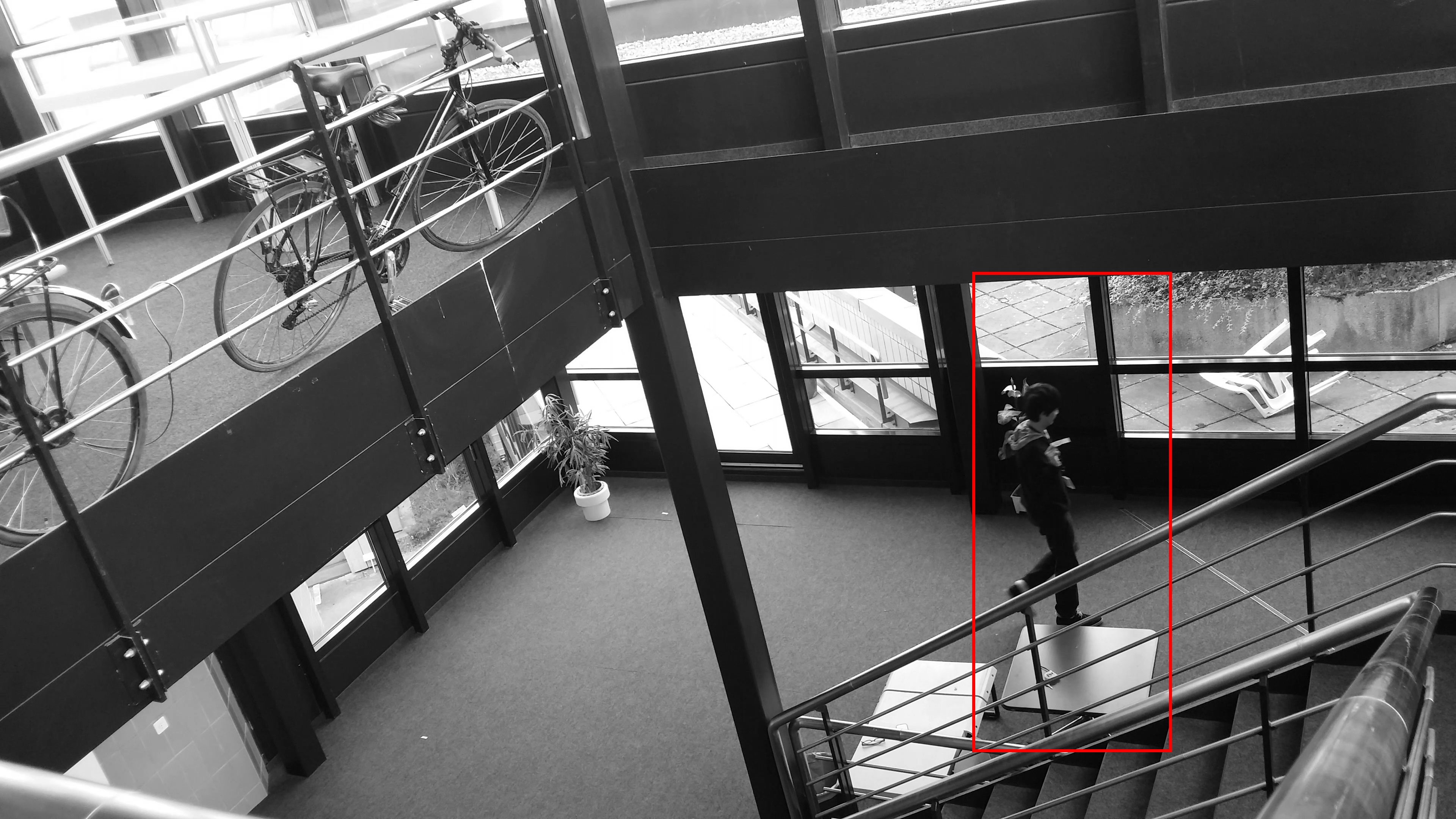} 
\put(48,-5){(a)}
\put(110,-5){(b)}
\put(135,-5){(c)}
\put(159,-5){(d)}
\end{overpic}
\includegraphics[height=.275\textwidth]{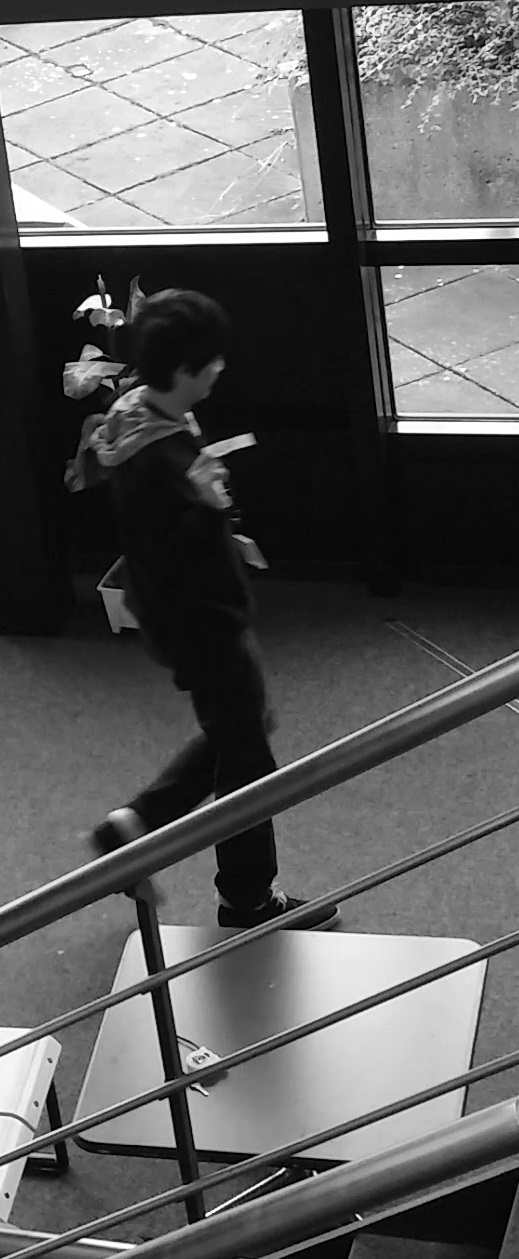} \hspace{-.05in}
\includegraphics[height=.275\textwidth]{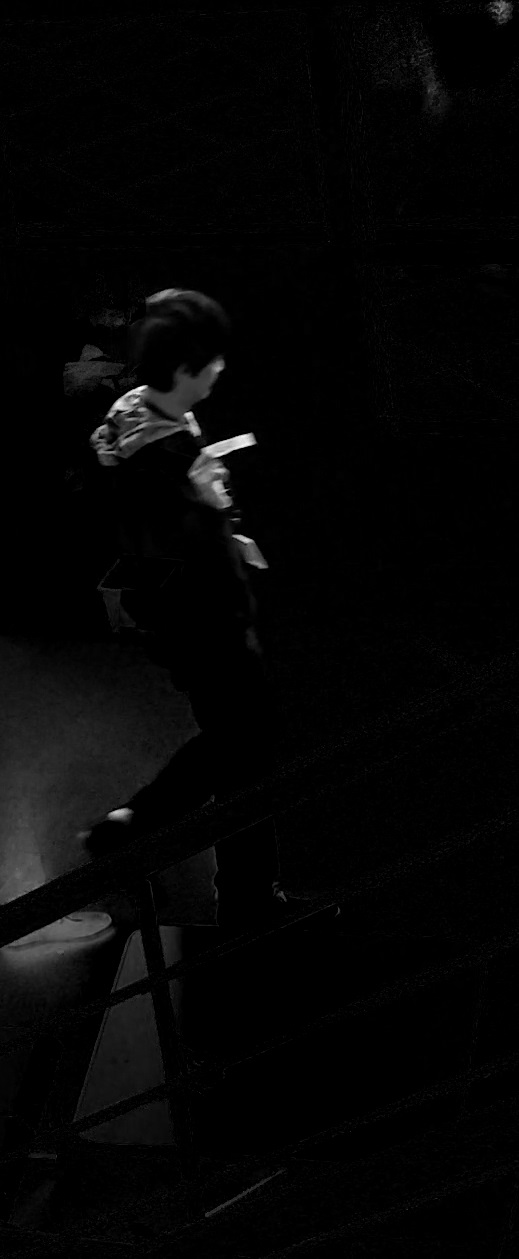}  \hspace{-.05in}
\includegraphics[height=.275\textwidth]{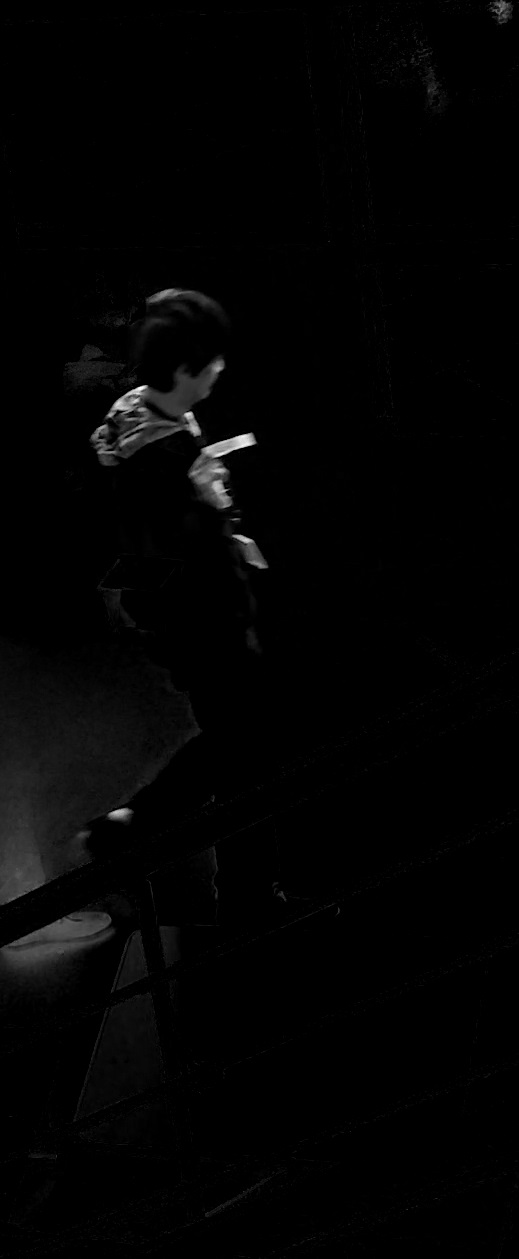} 
\end{center}
    \caption{Comparison of background subtraction on PEViD ``Walking Day Indoor 4"~\cite{KorUHD2014}, cropped.
(a) is the original frame, with a red rectangle indicating the boundaries of the non-zero region of the thresholded foreground.
(b) is a close-up of the original frame.
(c) is the foreground from DMD background subtraction on the CPU.
(d) is the foreground using our Streaming GPU DMD.
}
    \label{fig:pevid} 
\end{figure}

It is important to verify that the significant speed-up of the streaming GPU implementations does not come with an unacceptable loss in accuracy.  
Table \ref{tab:error} shows comparisons made between our streaming SVD and DMD output against Python implementations of the standard algorithms.
In both cases, the relative error is quite small, even for the largest input sizes.
The DMD comparison was made on the product of column of $\bPhi$ corresponding to the smallest absolute value in $\bLambda$.
The relative error is somewhat larger than that of the SVD, in part because of accumulated floating-point math errors from the GPU (e.g. \textit{fused multiply add} instructions).  
However, in many applications of DMD, such as video background modeling, this constitutes an acceptable error for the considerable speed-up, as downstream processing algorithms do not require machine precision. 

To show that the streaming GPU background subtraction is sufficiently accurate, we benchmark on the ``Video 003" from the Background Modeling Challenge dataset \cite{Vacavant2013}. The results are shown in Table \ref{tab:video}, and are consistent with other versions of the DMD\cite{erichson2015compressed}.  
We also provide a subjective comparison of background subtraction between our streaming GPU and the standard algorithm in Figure \ref{fig:pevid}.
The bottom row of Figure \ref{fig:pevid} shows, from left to right, a close-up of the original frame, the CPU foreground, the streaming GPU foreground and the difference between the two foregrounds.
The difference is small, and has little impact on the thresholded results.  

\section{Discussion and Outlook} \label{sec:discussion}
We have developed and analyzed streaming singular value and dynamic mode decomposition algorithms and their GPU implementations.
In addition, we show performance benefits for streaming video background subtraction.
In all cases, a large number of calculations are able to be carried forward from frame to frame by exploiting the structure of the method of snapshots SVD.
This allows both the SVD and DMD to process large data streams in real-time, whether for video or otherwise.
We have evaluated the proposed algorithms on multiple datasets, demonstrating the significantly improved computational performance for stream processing with negligible loss in accuracy.  
Our C++ and CUDA implementation will be made available under an open-source license.

The results of our performance comparison suggest that streaming algorithms are favorable, regardless of whether a GPU is available on a target platform. 
Additionally, significant speed-ups are possible at smaller data sizes once faster transfers are available to and from a GPU.
While not suitable for extreme-precision applications, we believe our streaming SVD and DMD algorithms provide a valuable improvement for many applications due to their improved computational performance.  
The small loss in accuracy was shown to be negligible for video background modeling applications. 

There are a number of interesting future directions that may arise from this work.  
One could modify the streaming algorithms shown here to support dynamic updating with more than one column at a time; when data inputs slow down, the number of new columns processed may be increased to catch up, and vice versa.
This dynamic streaming update could help to recover from a build-up of columns waiting to be processed for a long-running instance of the streaming SVD or DMD.
A streaming input build-up could also be used instead of waiting for enough initial inputs for the first SVD or DMD.
This would instead pre-allocate the maximum matrix size, but start the algorithm with only 2 columns.
Until the matrix is filled, the new columns would be appended without erasing the oldest.
In using the streaming DMD for background subtraction, the algorithm could be modified to use some small subset of background DMD modes rather than just the single slowest changing mode, as suggested in Erichson et al.~\cite{erichson2015compressed}.
This would incur a performance penalty, but could improve results as well.
Lastly, our method could be joined with other modified DMD algorithms, such as compressed or multi-resolution in order to improve performance.

The emergence of the big data era across the physical, biological, social and engineering sciences has severely challenged our ability to extract meaningful features from data in a real-time manner.   
Critical technologies such as LIDAR, 4K video streams, computer vision, high-fidelity numerical simulations, sensor networks, brain-machine interfaces, internet of things, and augmented reality will all depend on scalable algorithms that can produce meaningful decompositions of data in real time.
Failure to compute data streams in real time results in a data mortgage~\cite{Polagye2014DOE} whereby the cost of collection and storage limits the available resources to analyze and extract features.   
We are already seeing this across the sciences where massive data sets are collected and stored, yet remained un-mined for informative features and/or critical information for automated decision making processes.   
The streaming technique presented here provides a mathematical architecture for real-time processing of data and extraction of features.  The method is adaptive, efficient, parallelizable and scalable, potentially enabling a host of applications currently beyond the capabilities of standard techniques.

\section*{Acknowledgements}

We gratefully acknowledge valuable discussions with Zhe Bai and Ben Erichson.
SLB and SDP acknowledge support from the Department of Energy under award DE-EE0006785 and thank NVIDIA for providing a Tesla K40 GPU for this research.  
SLB also acknowledges support from the Boeing corporation and SDP acknowledges support from the Mary Gates Scholarship.  

\bibliographystyle{acm}
\bibliography{references}

\end{document}